\input amstex


\def\b1{\text{\bf 1}}

\def\BC{{\Bbb C}}

\def\BP{{\Bbb P}}
\def\BZ{{\Bbb Z}}

\def\CB{{\Cal B}}
\def\CD{{\Cal D}}
\def\CDiff{{\Cal{D}}iff}
\def\CE{{\Cal E}}

\def\CHom{{\Cal{H}}om}

\def\CM{{\Cal M}}

\def\CO{{\Cal O}}

\def\CU{{\Cal U}}
\def\CV{{\Cal V}}

\def\fg{{\frak g}}

\def\hfg{\widehat{\frak g}}

\def\Res{\text{Res}}

\def\tx{\tilde x}

\def\ty{\tilde y}


\def\btu{\bigtriangleup}

\def\iso{\buildrel\sim\over\longrightarrow} 

\def\lra{\longrightarrow}

\parskip=6pt

\documentstyle{amsppt}
\document
\magnification=1200
\NoBlackBoxes


\centerline{\bf Chiral Poincar\'e duality}

\bigskip
\centerline{Fyodor Malikov, Vadim Schechtman}
\bigskip

\bigskip

\bigskip

{\bf 1.} Let $X$ be a smooth algebraic variety over $\BC$. In the note
[MSV] we introduced a sheaf of vertex superalgebras $\Omega_X^{ch}$ 
on $X$. (Below we will often omit the prefix "super"; we will live mainly 
in the $\BZ/2$-graded world, the tilde over a letter will denote its 
parity.)  
This sheaf has a $\BZ\times\BZ_{\geq 0}$-grading 
$$
\Omega^{ch}_X=\oplus_{p\in\BZ,\ i\in\BZ_{\geq 0}}\ \Omega^{ch,p}_i
$$ 
by {\it fermionic number} $p$ and {\it conformal weight} $i$. 
The $\BZ/2$-grading is $p\ mod(2)$.  
The conformal weight zero part $\Omega^{ch}_0=\oplus_p\ \Omega^{ch,p}_0$ 
is identified 
with the usual de Rham algebra $\Omega_X=\oplus_p\ \Omega^p_X$ 
of differential forms. The sheaf $\Omega^{ch}_X$ will be called 
the {\it chiral de Rham algebra} of $X$. 
 
Each component $\Omega^{ch,p}_i$ carries a canonical finite 
filtration whose factors are locally free $\CO_X$-modules of finite rank. 
However, there is {\it no} natural $\CO_X$-module structure  
on $\Omega^{ch}_X$ itself. Let us consider the cohomology 
$$
H^*(X,\Omega^{ch}_X)=\oplus_{p, q, i}\ H^q(X,\Omega^{ch,p}_i)
$$
This is a conformal vertex superalgebra. The $\BZ/2$-grading 
is $(p+q)\ mod(2)$.  
When $X$ is 
complete we call this algebra   
the {\it chiral Hodge cohomology} of $X$. The conformal weight 
zero part of it coincides with the usual 
Hodge cohomology algebra $H^*(X,\Omega_X)=\oplus_{p,q}\ H^q(X,\Omega^p_X)$. 
If $X$ is Calabi-Yau then $H^*(X,\Omega^{ch}_X)$ is a $N=2$ superconformal 
vertex algebra, in the sense of [K], 5.8. 

We have the K\"unneth formula: for any two 
smooth varieties $X, Y$ a canonical isomorphism 
of conformal vertex superalgebras
$$
H^*(X\times Y,\Omega^{ch}_{X\times Y})=H^*(X,\Omega^{ch}_X)
\otimes H^*(Y,\Omega^{ch}_Y)
\eqno{(1.0)}
$$  

From now on we assume that {\it $X$ is complete}, unless specified otherwise. 
The "chiral Hodge --- de Rham spectral sequence" degenerates not 
at $E_1$ but at $E_2$. Namely, 
the chiral de Rham differential $d^{ch}_{DR}$ on $\Omega^{ch}_X$ induces 
a differential 
$$
Q:\ H^q(X,\Omega^{ch,p}_i)\lra H^q(X,\Omega^{ch,p+1}_i)
$$
and the cohomology of $H^*(X,\Omega^{ch})$ with respect 
to $Q$ is equal to $H^*(X,\Omega_X)$. Indeed, as in the proof 
of [MSV], Theorem 2.4, the operator $G_0$ is a zero homotopy   
on the components of nonzero conformal weight. 
The "chiral de Rham cohomology" coincides with the usual 
de Rham cohomology.

    We can consider the similar sheaf $\Omega^{ch,an}_X$ over the 
corresponding analytic variety $X^{an}$. 
We have canonical isomorphism
$$
H^*(X,\Omega^{ch}_X)\iso H^*(X^{an},\Omega^{ch,an}_X)
$$
Indeed, we have an obvious map from the left hand side to the right 
hand side, compatible with the above mentioned filtrations, 
and it is an isomorphism by GAGA and five-lemma.  

The aim of this note is to prove 

{\bf 1.1. Theorem.} {\it 
The space $H^*(X,\Omega^{ch}_X)$ carries   
a canonical non-degenerate symmetric bilinear form
$$
\langle\ ,\ \rangle:\ H^*(X,\Omega^{ch}_X)\times H^*(X,\Omega^{ch}_X)
\lra\BC
\eqno{(1.1)}
$$ 
This form makes the components of different conformal 
weights orthogonal and identifies   
$$
H^q(X,\Omega^{ch,p}_i)^*=H^{n-q}(X,\Omega^{ch,n-p}_i)
\eqno{(1.2)}
$$
Here $n=dim(X)$.  
The restriction of this form to the conformal weight zero 
component coincides with the usual Poincar\'e pairing.}

The symmetry of (1.1) is understood in the $\BZ/2$-graded sense.

{\bf 1.2. Remark.} Recall that we have defined in [MS], Part II, 1.4 
for every conformal vertex algebra $V$ a canonical Lie algebra 
antiinvolution 
$$
\eta:\ Lie(V)\lra Lie(V)
\eqno{(1.3)}
$$
Here $Lie(V)$ denotes the Lie algebra of Fourier components 
of the fields of $V$. 

It follows from the construction of (1.1) that this pairing is 
$\eta$-contravariant, i.e. 
$$
\langle x_{(n)}y,z\rangle = (-1)^{\tx\ty}\langle y,\eta(x_{(n)})z\rangle
\eqno{(1.4)}
$$  

Of course the Poincar\'e duality in the usual Hodge cohomology  
is an immediate consequence of the Serre duality 
$H^i(X,E)^*=H^{n-i}(X,E^o)$ where for a vector  
bundle   
$E$ on $X$, $E^o$ denotes the dual bundle $\CHom_{\CO_X}(E,\omega_X)$,
$\omega_X:=\Omega^n_X$. Similarly, we deduce Theorem 1.1  
from  the corresponding local statement, 
see Theorem 8.1. To formulate
it, we need to define the dual of the sheaf $\Omega^{ch}_X$; 
this is not immediate since $\Omega^{ch}_X$ is not an $\CO_X$-module. 
We define the dual using M.~Saito's 
language of induced $D$-modules, [S].

In fact we do more:  we introduce
a suitable category of "restricted"
(in the sense of [MSV])  $\Omega^{ch}_{X}$-modules along with a duality
functor on it and prove
 a  general statement, which can be thought of as a chiral analogue of  
 Serre duality,
see no. 11, formula (11.6).

Both the chiral Serre duality and Theorem 8.1  are
 consequences of  Theorem 10.1, which may be of independent 
interest. It says that the "Weyl module" functor is an equivalence 
between the categories of $\CD_{\Omega_X}$-modules and restricted
 $\Omega^{ch}_X$-modules. Theorem 11.2 adds that this equivalence
preserves the duality functor.

In no. 12 we present 
an alternative proof of 1.1 for $X=\BP^1$ and in no. 13 
derive some consequences about the structure 
of $H^*(\BP^1,\Omega^{ch}_{\BP^1})$.  

{\bf 2.} First let us recall Saito's theory. By  
a $D$-module on $X$ we mean a right $\CD_X$-module 
quasicoherent over $\CO_X$. The category of $D$-modules 
on $X$ will be denoted $\CM(X)$.  
For a $D$-module $M$, let $DR(M)$ denote 
its de Rham complex 
$$
DR(M):\ 0\lra M\otimes_{\CO_X}\Lambda^n\Theta_X\lra
\cdots\lra M\otimes_{\CO_X}\Theta_X\lra M\lra 0
\eqno{(2.1)}
$$
(we regard it as sitting in degrees $-n,\ldots,0$). Here 
$\Theta_X$ is the tangent sheaf. The de Rham cohomology 
$H^*_{DR}(X,M)$ is defined as the hypercohomology 
$H^*(X,DR(M))$.    
Set  
$$
h(M)=H^0(DR(M))=M/M\Theta_X
\eqno{(2.2)}
$$
For a quasicoherent $\CO_X$-module $P$, set $P^\sim:=P\otimes_{\CO}\CD_X$; 
this is a 
$D$-module (the action of $\CD_X$ is induced by the right 
action of $\CD_X$ on itself). A $D$-module isomorphic to 
$P^\sim$ for some $P$ is called 
{\it induced}. 

The de Rham complex 
$DR(P^\sim)$ is a left resolution of $P$; more precisely, 
we have a canonical arrow
$$
\nu_P:\ DR(P^\sim)\lra P
\eqno{(2.3)}
$$
sending $p\otimes\partial\in P\otimes\CD_X=DR^0(P^\sim)$ to 
$p\partial$, and $\nu_P$ is a quasiisomorphism. In particular, 
$h(P^\sim)=P$. As a consequence, we have a canonical isomorphism  
$$
H^*_{DR}(X,P^\sim)=H^*(X,P)
\eqno{(2.4)}
$$
A morphism of $D$-modules $f:\ P^\sim\lra Q^\sim$ induces  
a morphism of sheaves $h(f):\ P\lra Q$.  
One checks that $h(f)$ is a 
differential operator 
and this gives an isomorphism    
$$
Hom_{\CD_X}(P^\sim,Q^\sim)=Diff(P,Q)
\eqno{(2.5)}
$$
where $Diff(P,Q)$ denotes the space of differential operators, 
in the sense of Grothendieck, acting from $P$ 
to $Q$, cf. [S], 1.20.   

{\bf 3. Definition.} {\it A $D$-bundle on $X$ is a locally free right $\CD_X$-
module of finite rank.} 

The $D$-bundles form a full subcategory $\CD-\CB un(X)$ of $\CM(X)$.

{\bf 4.} Let $P$ be a sheaf of $\BC$-vector spaces on $X$  
We will call $P$ an {\it differential bundle} if it satisfies the 
propertiy (Diff) below. 

First let us formulate a weaker property   

(S) There exists a  a
Zariski open covering  
$\CU=\{ U\}$ of $X$ and $\BC$-linear isomorphisms of sheaves 
$$
s_U:\ P_U\iso E^U, 
\eqno{(4.1)}
$$
$U\in\CU$, for some vector bundles ($:=$ locally free $\CO_X$-modules  
of finite rank) $E^U$ over $U$  
Here the subscript $_U$ denotes the restriction to $U$. 

Let us call a collection of isomorphisms (4.1) a local trivialization 
of $P$. On   
the pairwise intersections, we get the isomorphisms 
$$
c_{UV}:=s_Vs_U^{-1}:\ E^U_V\iso E^V_U
\eqno{(4.2)}
$$ 
satisfying an obvious cocycle condition. 
Now we formulate the property (Diff) which strengthens (S):  

(Diff) There exists a local trivialization $\{ s_U\}$ such that 
the corresponding 
\newline $c_{UV}\in Diff(E^U_V,E^V_U)$.

{\bf 5.} If $M$ is a $D$-bundle then $h(M)$ is obviously
a differential bundle. 

Conversely, given an differential bundle $P$, choose a local 
trivializtion satisfying (Diff). Over each $U\in\CU$ one can 
from the induced module $E^{U\sim}$. Let us glue 
them together using the Cech $1$-cocycle $c^\sim=(c_{UV}^\sim)$ where
$$
c_{UV}^{\sim}\in Hom_{\CM(U\cap V)}(E^{U\sim}_V,E^{V\sim}_U)
$$ 
corresponds to $c_{UV}$ via (2.5). We get a  $D$-bundle  
$P^\sim$. 

For two differential bundles $P,P'$, we define the space 
$Diff(P,P')$ of differential operators as 
$$
Diff(P,P')=Hom_{\CM (X)}(P^\sim,P^{\prime\sim})
$$
Note that $Diff(P,P')$ is canonically a subspace of $Hom_\BC(P,P')$. 
This way we get a category $\CD iffbun(X)$ of differential bundles.  

The correspondences $M\mapsto h(M),\ P\mapsto P^\sim$ give rise 
to the quasiinverse equivalences between $\CD-\CB un(X)$ and $\CD iffbun(X)$.   

The obvious morphism $\nu:\ DR(P^\sim)\lra P$ is a 
quasiisomorphism for each differential bundle $P$. 
Consequently, we have canonically
$$
H^*_{DR}(X,P^\sim)=H^*(X,P)
\eqno{(5.1)}
$$

{\bf 6.} For a finite non-empty set $I$, an $I$-family 
$\{ P_i\}$ of differential bundles and a differential bundle $Q$ let 
$Diff_I(\{ P_i\},Q)$ denote the subspace of the space of 
maps $\otimes_{I\BC} P_i\lra Q$ 
which are differential operators by each argument (when 
all but one arguments are fixed). 
These spaces define a 
pseudo tensor structure on $\CD iffbun(X)$, in the sense of 
Beilinson-Drinfeld, [BD]. 

On the other hand, the category $\CD-\CB un(X)$ carries 
a pseudo tensor structure induced from the 
$*$-pseudo tensor structure on $\CM(X)$  
introduced in {\it op. cit.}, 2.2.3. 
More precisely, for an $I$-family of  
$D$-bundles $\{ L_i\}$ and a $D$-bundle $M$, 
we define the space of polylinear operations 
$P^*_I(\{ L_i\},M)$ as the space of all $D$-module maps  
$$
\boxtimes_I:\ L_i\lra \Delta^{(I)}_*M
$$
Here $\Delta^{(I)}:\ X\lra X^{(I)}$ 
is the diagonal embedding. 

The functor $h$ identifies both pseudo tensor 
structures: we have canonically
$$
Diff_I(\{ P_i\},Q)=P^*_I(\{ P^\sim_i\},Q^\sim)
\eqno{(6.1)}
$$

{\bf 7. Duality.} Recall the duality functor for 
$D$-modules. Consider the sheaf  
$\omega_X^\sim:=\omega_X\otimes_{\CO_X}\CD_X$. It carries two commuting 
structures of a right $\CD_X$-module: the first one 
coming from the tensor product of a right and a left 
$\CD_X$-module, and the second one appearing from the 
right $\CD_X$-module structure on $\CD_X$. Note that 
according to a lemma by Saito, [S], Lemma 1.7, there is a canonical 
involution on the above sheaf which interchages  
two $\CD_X$-module structures. 

Let $M$ be a $D$-bundle. Set   
$$
M^o:=\CHom_{\CD_X}(M,\omega_X^\sim)
\eqno{(7.1)}
$$
where $\CHom$ is taken with respect to the first 
$\CD_X$-module structure on $\omega_X^\sim$, and 
the right $\CD_X$-action on it is induced by the second 
structure; it is also a $D$-bundle.

{\bf 7.1. Lemma.} {\it One has canonical isomorphisms 
$H^i_{DR}(X,M)^*=H^{n-i}_{DR}(X,M^o)$ where $n=dim(X)$.}   

This is the duality theorem for $D$-modules, cf. [S].  

Note that if $M=E^\sim$ where $E$ is a vector bundle 
then $M^o=E^{o\sim}$ where $E^o=\CHom_{\CO_X}(E,\omega_X)$.

{\bf 7.2. Definition.} {\it Let $P$ be a differential bundle. 
The dual  differential bundle $P^o$ is defined by 
$P^o=h(P^{\sim o})$.} 

One can give a more direct definition of $P^o$ using the gluing 
functions. Namely, choose a local trivialization $\{ s_U\}$ as in 4.1 (S) 
satisfying (Diff),  
with the corresponding cocycle $c=(c_{UV})$. Recall that for 
two vector bundles $E, F$ we have a canonical isomorphism 
$$
Diff(E,F)=Diff(F^o,E^o)
\eqno{(7.2)}
$$
For example, if $E=F=\CO_X$ then (7.2) amounts to the usual 
correspondence between left and right $\CD_X$-modules.

 Now, the dual differential bundle $P^o$ is glued by means of the 
dual Cech cocycle $c^o=(c_{UV}^o)$ where 
$c_{UV}^{o}\in Diff(E^{Vo}_U,E^{Uo}_V)$ 
corresponds to $c_{UV}$ via (7.2). 

Note that we have canonically 
$$
E^{oo}=E
\eqno{(7.3)}
$$   
 
Lemma 7.1 along with (5.1) implies 

{\bf 7.3. Lemma} (Serre duality) {\it One has canonical isomorphisms
 $H^i(X,P)^*=H^{n-i}(X,P^o)$}.

{\bf 8.} The arguments of 
[MSV], 6.10 (cf. also [MS], Part I) show that a choice of \'etale coordinates 
on a Zariski open $U\subset X$ gives a trivialization  
of the sheaves $\Omega^{ch,p}_{iU}$, and one sees that the transition 
functions are differential operators. 
It follows that the sheaves $\Omega^{ch,p}_i$ carry a 
canonical structure of differential bundles.  

{\bf 8.1. Theorem.} {\it For all $p,i$ there exist canonical isomorphisms 
$$
\chi_i^p:\ (\Omega^{ch,p}_i)^o\iso \Omega^{ch,n-p}_i
\eqno{(8.1)}
$$ 
For $i=0$ the isomorphisms} (8.1) {\it are induced by the wedge 
products of differential forms.

These isomorphisms are symmetric in the following sense: 
$$
\chi_i^{p o}=\chi_i^{n-p}
\eqno{(8.2)}
$$}

The proof is given after the proof of Theorem 11.2 below. 

  Theorem 8.1 and Lemma 7.2 3 immediately imply Theorem 1.1. 

{\bf 9.} Note that the sheaf $\Omega^{ch}_X$ is a 
vertex algebra in the Diff-pseudo tensor 
structure defined in 6 (a {\it vertex$^{Diff}$ algebra}). 
This means that the operations 
$$
_{(n)}:\ \Omega^{ch}_X\times\Omega^{ch}_X\lra \Omega^{ch}_X
$$
belong to $Diff_{\{1,2\} }(\{\Omega_X^{ch},\Omega^{ch}_X\},
\Omega^{ch}_X)$. 

Equivalently, the $D$-module $\Omega^{ch\sim}_X$ is 
a {\it vertex$^*$ algebra} (a vertex algebra 
in the $*-$ pseudo tensor structure).

 {\bf 10.} Let us consider the de Rham algebra of differential 
 forms $\Omega_X=\oplus\ \Omega_X^p$. Let $\CD_{\Omega_X}$ be 
 the superalgebra  
 of differential operators on $\Omega_X$. Let $\CD_{\Omega_X}-Mod$ 
 denote the category of left $\CD_{\Omega_X}$-modules (everything here  
 is $\BZ/2$-graded). Let 
 $\Omega^{ch}_X-Mod$ denote the category of restricted 
 $\Omega_X^{ch}$-modules ("restricted" means that the modules are 
 graded by conformal weight and there are no components of 
 negative weight).  
  
 Recall that we defined in [MSV], 6.11, [MS], I.4.5 the {\it Weyl} functor
 $$
 W_\Omega:\ \CD_{\Omega_X}-Mod\lra \Omega^{ch}_X-Mod
 \eqno{(10.1)}
 $$
 This functor is simply the left adjoint to the functor of taking 
 the conformal weight zero component. 
 More explicitely, let $U\subset X$ 
 be a sufficiently small Zariski open. 
 A choice of \'etale coordinates on $U$ provides us with two things: 
it makes   $\Gamma(U,\Omega^{ch}_{X})$
into a supercommutative algebra, and it determines an embedding of
algebras $\Omega_{X}\hookrightarrow\Omega^{ch}_{X}$. We then set
$$
\Gamma(U,W_{\Omega}(V))=
\Gamma(U,\Omega^{ch}_{X})\otimes_{\Gamma(U,\Omega_{X})}V
\eqno{(10.2)}
$$
$\Gamma(U,W_{\Omega}(V))$ defined in this way is obviously a
$\Gamma(U,\Omega^{ch}_{X})$-module. One further checks that this
$\Gamma(U,\Omega^{ch}_{X})$-module structure is in fact independent of the
choice of coordinates and nicely agrees with localization. The
sheaf $W_{\Omega}(V)$ is then defined by picking a suitable affine atlas of $X$.

For example, $W_{\Omega}(\Omega_X)=\Omega_X^{ch}$. 
In general, we refer  
to  $W_{\Omega}(V)$ as a Weyl module over $\Omega^{ch}_{X}$, or simply
a chiral Weyl module.
 
 {\bf 10.1. Theorem.} {\it The functor $W_\Omega$ is an equivalence 
 of categories.}
 
 {\bf Proof.} By definition, if $\CV\in\Omega^{ch}_{X}-Mod$, then its
conformal weight 0 component, $\CV_{0}$, is a $\CD_{\Omega_{X}}$-module and
$W_{\Omega}(V)_{0}=V$. Therefore,
it suffices to show that for any  $\CV\in\Omega^{ch}_{X}-Mod$
there is an isomorphism $W_{\Omega}(\CV_{0})\rightarrow \CV$.

To construct a map $W_{\Omega}(\CV_{0})\rightarrow \CV$ observe that 
by definition $W_{\Omega}(\CV_{0})$ has the following universality 
property:

for any   $\CU\in\Omega_{X}^{ch}-Mod$, any $V\in \CD_{\Omega_{X}}-Mod$,
and any morphism of
 $\CD_{\Omega_{X}}$-modules
$$
f: V\rightarrow \CU_{0}
$$
there is a unique morphism
 $$
\hat{f}: W_{\Omega}(V)\rightarrow\CU
$$
such that the restriction of $\hat{f}$ to $V\subset W_{\Omega}(V)$
equals $f$.  Therefore we get the map
$$
\widehat{id}: W_{\Omega}(\CV_{0})\rightarrow \CV.
$$

To prove injectivity and surjectivity of $\widehat{id}$ we introduce, for
any $\CV\in\Omega^{ch}_{X}-Mod$, the subsheaf of {\it singular vectors},
$Sing(\CV)$, to be defined as follows:
$$
\Gamma(U, Sing(\CV))=\{v\in \Gamma(U, \CV): x_{i}v=0\text{ for all }
x\in\Gamma(U, \Omega_{X}^{ch}), i>0\}.
\eqno{(10.3)}
$$

{\bf 10.2. Lemma.}  

(i) {\it $Sing(\CV)\neq 0$ for any $\CV\in\Omega_{X}^{ch}-Mod$.}

(ii) $Sing(\CV)\subset \CV_0$.

    This lemma allows us to complete the proof of the theorem at once.
By (ii) $Sing(Ker(\widehat{id}))=Ker(\widehat{id})_{0}$, which
equals $Ker(id)$ and is, therefore, 0.
 Hence, by (i), $Ker(\widehat{id})=0$.
Replacing $Ker$ with $Coker$ in this argument, we get that
$Coker(\widehat{id})=0$

Let us finally prove the lemma. Item (i) is an obvious consequence of
the restrictedness property: it is enough to observe that
$\CV_{i}\subset Sing(\CV)$ if $i\geq 0$ is the minimum number satisfying
$V_{i}\neq 0$.

 As to (ii), we remind the reader that
 conformal weights  are  eigenvalues of one of the Virasoro generators,
$L_{0}$, which is given locally (over, say, a formal polydisk) by the formula

$$
L_{0}=\sum_{i,k}i:a_{i}^{k}b_{-i}^{k}:+i:\phi_{-i}^{k}\psi_{i}^{k}:.
$$
Comparing with (10.3) we see that because of the coefficient $i$,
$L_{0}$ acts as 0 on $Sing(\CV)$.

 {\bf 11. Chiral Serre duality.} We define a {\it $\Omega^{ch}_X$-bundle} to be
 an $\Omega^{ch}_X$-module $\CE$ such that it is a differential 
 bundle and all operations
 $$
 _{(n)}:\ \Omega^{ch}_X\otimes \CE\lra \CE
 $$
 belong to $Diff_{\{1,2\}}(\{\Omega^{ch}_X,\CE\},\CE)$. 
 
 Let $\CE$ be a $\Omega^{ch}_X$-bundle. 
 Consider the
{\it restricted} dual differential 
 bundle $\CE^o$. By this we mean the following:
being graded by conformal weight $\CE$ is a direct sum of
differential bundles $\CE=\oplus_{i\geq 0}\CE_{i}$ and we set
$$
\CE^{o}=\oplus_{i\geq o}\CE_{i}^{o}.
\eqno{(11.1)}
$$
(We do not change the notation hoping that this will not lead to a confusion:
it is a general principle that in the realm of modules over a vertex algebra
a dual means a restricted dual.)

 We want to introduce a canonical structure 
 of an $\Omega^{ch}_X$-bundle on it. 
 
 Note that if $M$ is a module over a vertex algebra $V$ then 
 $M$ is automatically a $Lie(V)$-module; the converse is not in general true.
 However, a $Lie(V)$-module structure on $M$ may come from 
 at most one $V$-module srtucture.

 Let us endow $\CE^o$ with a
 $Lie(\Omega^{ch}_X)$-module structure. First of all, by (2.5) the $D$-bundle
$\CE^{\sim o}= \CHom_{\CD_{X}}(\CE^{\sim},\omega^{\sim})$ equals
$\CDiff(\CE,\omega)$. Therefore it carries a canonical structure of a right
 $Lie(\Omega^{ch}_X)$-module defined by the formula
$$
(x_{n}f)(.)=(-1)^{\tilde{x}\tilde{f}}f(x_{n}.),
\eqno{(11.2)}
$$
and hence  a canonical structure of a left
 $Lie(\Omega^{ch}_X)$-module defined by the formula
$$
(x_{n}f)(.)=(-1)^{\tilde{x}\tilde{f}}f(\eta(x_{n}).),
\eqno{(11.3)}
$$
where $x$ is a local section of $\Omega^{ch}_{X}$. Here 
$$
\eta:\ Lie (\Omega^{ch}_X)\lra 
Lie(\Omega^{ch}_X)
$$
is the canonical antiinvolution (see 1.2). 

Second of all, this  $Lie(\Omega^{ch}_X)$-module structure descends to
the quotient $\CE^{o}=\CDiff(\CE,\omega)/\CDiff(\CE,\omega)\Theta_X$; 
this is because
the action of the tangent sheaf $\Theta_X$ commutes with the action of
 $Lie(\Omega^{ch}_X)$: $\Theta_X$ acts on the value of the function $f(.)$,
while $Lie(\Omega^{ch}_X)$ acts on its argument, see (11.2).

 {\bf 11.1. Claim.} {\it The above $Lie(\Omega^{ch}_X)$-module 
 structure on $\CE^o$ comes from the $\Omega^{ch}_X$-module 
 structure.}
 
  {\bf Proof.} Let temporarily $V$ be a vertex algebra.
To prove the claim we have to understand what is it that
singles out $V$-modules from the class of $Lie(V)$-modules. A pair
$(E,\rho)$ is a $Lie(V)$-module if $E$ is a vector space and
$\rho: Lie(V)\rightarrow End(E)$ is a Lie (super)algebra morphism.
In particular, for any $x\in V$ we have a family of operators
$\rho(x_{n})\in End(E)$. For $(E,\rho)$ to be
a $V$-module the two additional conditions are to be satisfied:

(A) For any $x\in V$ and $e\in E$, $\rho(x_{n})=0$ for all
$n>>0$.

(B) The operators $\rho(x_{n})\in End(E)$ satisfy the Borcherds
identities.

It follows from [K] Proposition 4.8 that those Borcherds identities that do not
follow from the Lie algebra structure on $Lie(V)$ follow from the
the relations:
$$
(x_{-\Delta_{x}}y)_{n}=\sum_{i\in\BZ}:x_{i}y_{n-i}:,\; n\in\BZ,
$$
where as usual $\Delta_{x}$ is a number such that $x\in V_{\Delta_{x}}$.

Therefore (B) is equivalent to

($\text{B}_{0}$) $(E,\rho)$ is a $Lie(V)$-module and
$$
\rho((x_{-\Delta_{x}}y)_{n})=\sum_{i\in\BZ}:\rho(x_{i})\rho(y_{n-i}):.
\eqno{(11.4)}
$$
The last formula means that each product $:\rho(x_{i})\rho(y_{n-i}):$
is ordered in the standard way and applied to any $e\in E$ from the right
to the left.

In the same way one compares right $V$-modules and right $Lie(V)$-modules
and concludes that a right $V$-module is a pair $(E,\rho)$ as above
satisfying the following conditions

($\text{A}^{r}$) For any $x\in V$ and $e\in E$, $\rho(x_{n})=0$ for all
$n<<0$.

($\text{B}_{0}^{r}$) $(E,\rho)$ is a right $Lie(V)$-module and
$$
\rho((x_{-\Delta_{x}}y)_{n})=\sum_{i\in\BZ}:\rho(x_{i})\rho(y_{n-i}):,
\eqno{(11.5)}
$$
where  each  product $:\rho(x_{i})\rho(y_{n-i}):$
is again ordered in the standard way but applied to any $e\in E$ from the left
to the right.

Having reviewed this undoubtedly well-known material
we cast a glance at (11.2) and convince ourselves that (11.2)
indeed determines a right $\Omega^{ch}_{X}$-module structure on
$\CE^o$: the restrictedness guarantees the condition ($\text{A}^{r}$),
while (14.5) holds simply because ($\text{B}_{0}$) holds for $\CE$.

Finally approaching (11.3) we see that (A) holds because
 ($\text{A}^{r}$) is satisfied for the right module structure
determined by (11.2)  and the fact that $\eta$ changes
the conformal weight to the opposite one. As to ($\text{B}_{0}$), it
is implied by the following easily checked property of the antiinvolution
$\eta$:
$$
\eta(\sum_{i\in\BZ}:x_{i}y_{n-i}:)=
\sum_{i\in\BZ}\eta(:x_{i}y_{n-i}:),
$$
where the action of $\eta$ on each monomial is as follows:
$$
\eta(x_{s}y_{t})=\eta(y_{t})\eta(x_{s}),
 \eta(y_{s}x_{t})=\eta(x_{t})\eta(y_{s}).
$$

 This defines an $\Omega^{ch}_X$-bundle structure on $\CE^o$. 
 By Lemma 7.3 we have 
 $$
 H^i(X,\CE^o)=H^{n-i}(X,\CE)^*
 \eqno{(11.6)}
 $$
 
 The wedge product $\Omega_X\times\Omega_X\lra\omega_X$ induces an 
 isomorphism of $\CO_X$-modules
 $$
 \Omega_X\iso\Omega_X^o
 \eqno{(11.7)}
 $$
 There is a unique left $\CD_{\Omega_X}$-module structure on $\Omega_X^o$ 
 such that (11.7) is an isomorphism of $\CD_{\Omega_X}$-modules. 
 
 More generally, for a left $\CD_{\Omega_X}$-module $E$, the 
 dual sheaf $E^o$ is canonically a right $\CD_{\Omega_X}$-module. 
 Indeed, we have an algebra homomorphism $\CD_{\Omega_X}\lra \CDiff(E)$,  
 hence $\CD_{\Omega_X}^o\lra \CDiff(E)^o=\CDiff(E^o)$ (cf. (7.2)). 
 On the other hand, the isomorphism 
 (11.7) induces an antiautomorphism
 $$
 \CD_{\Omega_X}\iso\CD^o_{\Omega_X}
 \eqno{(11.8)}
 $$
 Therefore, $E^o$ gets a canonical structure of a left 
 $\CD_{\Omega_X}$-module. This is the conformal 
 weight zero part of the definition of duality at the beginning 
 of this no. 
 
 {\bf 11.2. Theorem.} {\it The functor $W_\Omega$ commutes with 
 duality, i.e.  we have natural isomorphisms of $\Omega^{ch}_X$-modules
 $$
 W_\Omega(E^o)\iso W_\Omega(E)^o
 \eqno{(11.9)}
 $$}

{\bf Proof.}  By construction  $ W_\Omega(E)^o_{0}=E^{o}$. On the other
hand, due to (10.2) $ W_\Omega(E^o)_{0}=E^o$. Hence
$ W_\Omega(E^o)_{0}=  W_\Omega(E)^o_{0}$. But what Theorem 10.1 tells us
is that an $\Omega^{ch}_{X}$-module is uniquely determined by its conformal
weight zero component.
Therefore (11.9) immediately follows from Theorem 10.1. $\btu$

The isomorphisms (8.1) are a particular case of (11.9) with $E=\Omega_X$. 
By 10.1, it suffices to check the symmetry (8.2) on the conformal 
weight zero level, where it is evident. This proves Theorem 8.1. 

Theorem 1.1 is an immediate consequence of Theorem 8.1 and Lemma 7.3. 

{\bf 12.} Let $X=\BP^1$. In this case (1.2) reduces to the following 
$$
H^0(\BP^1,\Omega^{ch,p}_{\BP^1})^*=H^{1}(\BP^1,\Omega^{ch,1-p}_{\BP^1})
\eqno{(12.1)}
$$
For the sake of a mistrustful reader we present here a direct proof
of (12.1).

{\bf 12.1.} First of all, we explicitly describe the space
 $\Gamma(\BC^{*},\Omega^{ch}_{\BP^1})$. Consider the  Lie (super)algebra
$\Gamma$ on the even generators $a_{i},b_{i},\;i\in\BZ$ odd generators
$\phi_{i},\psi_{i},\;i\in\BZ$ and relations:

$$[a_{i},b_{-i}]=[\psi_{i},\phi_{-i}]=1,
\eqno{(12.2)}$$
all other brackets being equal 0.

This algebra is $\BZ$-graded (by conformal weight)
$$\Gamma=\oplus_{i\in\BZ}\Gamma_{i},$$
so that $\Gamma_{i}$ is linearly spanned by $x_{i}$, where $x$ is $a$, $b$,
$\phi$, or $\psi$. There arise four subalgebras

$$\Gamma_{>}=\oplus_{i>0}\Gamma_{i},
\Gamma_{<}=\oplus_{i<0}\Gamma_{i},
\Gamma_{\geq}=\oplus_{i\geq 0}\Gamma_{i},\; \Gamma_{0},$$
and the decomposition
$$
\Gamma=\Gamma_{<}\oplus\Gamma_{0}\oplus\Gamma_{>}
\eqno{(12.3)}
$$

For any Lie (super)algebra $\fg$, denote by $U(\fg)$ its universal
enveloping algebra. We have the extension of Lie 
algebras
$$
0\rightarrow\Gamma_{>}\rightarrow\Gamma_{\geq}\rightarrow\Gamma_{0}
\rightarrow 0,\eqno{(12.4)}
$$

If we fix $b$, a coordinate on $\BC^{*}$, and identify $b_{0}$ with $b$,
$a_{0}$ with $d/db$, $\phi_{0}$ with $db$, and finally $\psi_{0}$ with
the odd vector field $d/d(db)$, then $U(\Gamma_{0})$ gets identified with
the algebra of differential operators acting on $\Gamma(\BC,\Omega_{\BP^1})=
\Gamma(\BC,\Omega_{\BP^1}^{0}\oplus \Omega_{\BP^1}^{1})$. Hence the latter
space, as well as the space  $\Gamma(\BC^{*},\Omega_{\BP^1})$, becomes
a $\Gamma_{0}$-module, and, by pull-back due to (12.4), a
$\Gamma_{\geq}$-module. The inspection of the relevant definitions
in [MSV] shows that
$$
 \Gamma(\BC^{*},\Omega^{ch}_{\BP^1})=\text{Ind}_{\Gamma_{\geq}}^{\Gamma}
\Gamma(\BC^{*},\Omega_{\BP^1}).
\eqno{(12.5)}
$$

{\bf 12.2.} There is a natural pairing
$$
<.,.>:\Gamma(\BC^{*},\Omega_{\BP^1})\otimes\Gamma(\BC^{*},\Omega_{\BP^1})
\rightarrow\BC,\; <\nu,\mu>= \Res_{b=0}\nu\wedge\mu.
\eqno{(12.6)}
$$
It enjoys the following `contravariance'  properties
$$
<b_{0}\nu,\mu>=<\nu,b_{0}\mu>,\;
<a_{0}\nu,\mu>=-<\nu,a_{0}\mu>.
\eqno{(12.7)}
$$
(Similar equalities hold for the odd elements $\phi_{0},\psi_{0}$.)

The pairing (12.6) induces the following map
$$
\Gamma(\BC^{*},\Omega_{\BP^1})\rightarrow \Gamma(\BC^{*},\Omega_{\BP^1})^{*}.
\eqno{(12.8)}
$$

Well-known in representation theory is the operation of taking a
contragredient module. If we have a Lie algebra $\fg$ with an antiinvolution
and a $\fg$-module $M$ graded by finite dimensional subspaces, then
the contragredient $\fg$-module, $M^c$, is defined as the restricted dual
of $M$, the action of $\fg$ being equal to the canonical right action
twisted by the antiinvolution. Apply this construction to $\Gamma_{0}$
operating on $\Gamma(\BC^{*},\Omega_{\BP^1})$. The two necessary
structures are as follows: the antiinvolution is defined by
$$
\eta:\Gamma_{0}\rightarrow\Gamma_{0},\;\eta(b_{0})=b_{0},
\eta(a_{0})=-a_{0}, \eta(\phi_{0})=\phi_{0},
\eta(\psi_{0})=-\psi_{0}
\eqno{(12.9)}
$$
and the grading on $\Gamma(\BC^{*},\Omega_{\BP^1})$ is determined
by the condition
$$
deg(b_{0})=1,\; deg(a_{0})=-1.
\eqno{(12.10)}
$$
In this way we get the contragredient module $\Gamma(\BC^{*},\Omega_{\BP^1})^c$.
It is obvious that the map (12.8) gives an isomorphism of
$\Gamma_{0}$-modules
$$
\Gamma(\BC^{*},\Omega_{\BP^1})\rightarrow \Gamma(\BC^{*},\Omega_{\BP^1})^{c}.
\eqno{(12.11)}
$$
The same construction applies to the $\Gamma$-module
$ \Gamma(\BC^{*},\Omega^{ch}_{\BP^1})$. We, first, define an antiinvolution
$$
\eta:\Gamma\rightarrow\Gamma,\;\eta(b_{i})=b_{-i},
\eta(a_{i})=-a_{-i},\eta(\phi_{i})=\phi_{-i},
\eta(\psi_{i})=-\psi_{-i}.
\eqno{(12.12)}
$$
As to the grading on $ \Gamma(\BC^{*},\Omega^{ch}_{\BP^1})$ , we notice that
for any $i$ the subspace of conformal weight $i$,
$ \Gamma(\BC^{*},\Omega^{ch}_{\BP^1})$, is infinite dimensional and we cure
this by setting
$$deg(x_{i})=i\text{ if }i\neq 0,
$$
$$
deg(x_{0})=0\text{ unless }x=b\text{ or }a
$$
$$deg(b_{0})=-deg(a_{0})=1.
$$
In this way we get the contragredient module
$\Gamma(\BC^{*},\Omega^{ch}_{\BP^1})^c$. By definition
 $ \Gamma(\BC^{*},\Omega^{ch}_{\BP^1})^c$ inherits the grading by
conformal weight:
$$
 \Gamma(\BC^{*},\Omega^{ch}_{\BP^1})^c=
\oplus_{i\geq 0} \Gamma(\BC^{*},\Omega^{ch}_{\BP^1})^{c}_{i}.
\eqno{(12.13)}
$$
Due to (12.5) $\Gamma(\BC^{*},\Omega^{ch}_{\BP^1})_{0}=
\Gamma(\BC^{*},\Omega_{\BP^1})$ and by (12.13)
 the map (12.11) is actually an isomorphism of $\Gamma_{\geq}$-modules
$$
\Gamma(\BC^{*},\Omega^{ch}_{\BP^1})_{0}\rightarrow
\Gamma(\BC^{*},\Omega^{ch}_{\BP^1})^{c}_{0}.
\eqno{(12.14)}
$$
The universality property of induced modules implies that
the map(12.14) uniquely extends to a morphism of $\Gamma$-modules
$$
\Gamma(\BC^{*},\Omega^{ch}_{\BP^1})\rightarrow
\Gamma(\BC^{*},\Omega^{ch}_{\BP^1})^{c}.
\eqno{(12.15)}
$$
The latter map gives rise to the pairing 
$$
<.,.>:\Gamma(\BC^{*},\Omega_{\BP^1}^{ch})\otimes\Gamma(\BC^{*},\Omega_{\BP^1}^{c
h})
\rightarrow\BC,
\eqno{(12.16)}
$$
which has the following  contravariance property ( cf. (12.7))
$$
<x\nu,\mu>=(-1)^{\tilde{x}\tilde{\nu}}<\nu,\eta(x)\mu>.
\eqno{(12.17)}
$$

It is easy to see that the contragredient form (12.16) is uniquely
determined by the property (12.17). Thus we could have  defined
this form by demanding that (12.17) be valid, but then a certain
argument proving existence would have  been required; we chose instead
to  construct the map (12.15) using the properties of induction.

A closer look at the process of calculating the form (12.16) by the
repeated application of (12.17) shows that is is commutative:
$$
<\nu,\mu>=(-1)^{\tilde{\nu}\tilde{\mu}}<\mu,\nu>.
\eqno{(12.18)}
$$

{\bf 12.3.} We have everything ready for the proof of (12.1).
The space $\Gamma(\BC^{*},\Omega^{ch}_{\BP^1})$ has two subspaces:
$\Gamma(\BC,\Omega^{ch}_{\BP^1})$ and
$\Gamma(\BP^1 -\{0\},\Omega^{ch}_{\BP^1})$. To somewhat simplify the
notation set $V=\Gamma(\BC^{*},\Omega^{ch}_{\BP^1})$,
$V_{0}=\Gamma(\BC,\Omega^{ch}_{\BP^1})$,
$V_{\infty}=\Gamma(\BP^1-\{0\},\Omega^{ch}_{\BP^1})$.
One easily checks that the map (12.15) is actually an isomorphism;
therefore the pairing (12.16) is non-degenerate and  we use it to
identify $V$ with $V^{*}$.
By construction $V_{0}$ equals
its annihilator, $\text{Ann}V_{0}$. Since $\Gamma_{0}$-invariance
allows us to interchange 0 and $\infty$, $V_{\infty}$ is also
equal to $\text{Ann}V_{\infty}$. We now compute in a rather standard
manner
$$
H^{0}(\BP^1,\Omega^{ch}_{\BP^1})^{*}=(V_{0}\cap V_{\infty})^{*}=
$$
$$
=V/(\text{Ann}V_{0}+\text{Ann}V_{\infty})
=V/(V_{0}+V_{\infty})=H^{1}(\BP^1, \Omega^{ch}_{\BP^1}).
$$ 

{\bf 13.} Let us deduce some corollaries from the previous construction. 
The form (12.16) has several attractive properties. We have already 
noted that it is symmetric (12.18) and contravariant (12.17). 

Let $\fg=sl(2)$. 
We have proven in [MS], Part III, \S 1 that the affine Lie algebra 
$\hfg$ acts canonically on the sheaf $\Omega^{ch}_{\BP^1}$.

{\bf 13.1. Claim.} {\it The pairing} (12.16) {\it is $\hfg$-contravariant.} 

In fact, (12.17) implies that (12.16) is contravariant with respect to
$Lie(\Gamma(\BC^{*},\Omega^{ch}_{\BP^1}))$; therefore it is also with
respect to $\hfg$ since the latter acts by means of an embedding
$\hfg\hookrightarrow Lie(\Gamma(\BC^{*},\Omega^{ch}_{\BP^1}))$.

 As we  noted in {\it op. cit.}, Part III, 2.1, the space 
 $H^0(\BP^1,\Omega^{ch}_{\BP^1})$ is the maximal $\fg$-integrable 
$\hfg$-submodule 
 of $H^0(\BC^*,\Omega^{ch}_{\BC^*})$. 
 
 {\bf 13.2. Corollary.} {\it The $\hfg$-module $H^1(\BP^1,\Omega^{ch}_{\BP^1})$
 is the maximal $\fg$-integrable quotient of $H^0(\BC^*,\Omega^{ch}_{\BC^*})$.}
 
 This fact, conjectured by B. Feigin, was the starting point 
 of the present note.

{\bf 13.3.} We presented in  {\it op. cit.}, Part III, 2.2 a rather explicit
description of  $H^0(\BP^1,\Omega^{ch}_{\BP^1})$ as a $\hfg$-module. Now (12.1)
and 13.1 provide us with no less explicit description of
 $H^1(\BP^1,\Omega^{ch}_{\BP^1})$. We leave the details for the interested
reader.


\bigskip
\centerline{\bf References}
\bigskip


[BD] A.~Beilinson, V.~Drinfeld, Chiral algebras I, Preprint.





[K] V.~Kac, Vertex algebras for beginners, University Lecture Series, 
{\bf 10}, American Mathematical Society, Providence, RI, 1997.

[MS] F.~Malikov, V.~Schechtman, Chiral de Rham complex. II, 
{\it D.B.~Fuchs' 60-th Anniversary volume} (1999), to appear; 
math.AG/9901065.  

[MSV] F.~Malikov, V.~Schechtman, A.~Vaintrob, Chiral de Rham complex,
{\it Comm. Math. Phys.} (1999), to appear. 

[S] M.~Saito, Induced $D$-modules and differential complexes, 
{\it Bull. Soc. Math. France}, {\bf 117} (1989), 361-387.

\bigskip

\bigskip

F.M.: Department of Mathematics, University of Southern California, 
Los Angeles, CA 90089, USA;\ fmalikov\@mathj.usc.edu 

V.S.: Department of Mathematics, Glasgow University, 
15 University Gardens, Glasgow G12 8QW, United Kingdom;\ 
vs\@maths.gla.ac.uk

\enddocument